\documentclass[11pt]{article}
\usepackage{latexsym,amssymb,amsmath,amsfonts,amsthm}

\makeatletter
\newif\ifmsbmloaded@
\def\loadmsbm{\msbmloaded@true
  \font\tenmsb=msbm10 scaled 1\@ptsize00
  \font\sevenmsb=msbm7 scaled 1\@ptsize00
  \font\fivemsb=msbm5 scaled 1\@ptsize00
  \alloc@8\fam\chardef\sixt@@n\msbfam
  \textfont\msbfam=\tenmsb
  \scriptfont\msbfam=\sevenmsb
  \scriptscriptfont\msbfam=\fivemsb
  }
\def\nonmatherr@#1{\errmessage%
{LateX error: \string#1\space allowed only in math mode}}
\def\Bbb{\relax\ifmmode\expandafter\Bbb@\else
  \expandafter\nonmatherr@\expandafter\Bbb\fi}
\def\Bbb@#1{{\Bbb@@{#1}}}
\def\Bbb@@#1{\fam\msbfam\relax#1}
\@addtoreset{equation}{section}

\makeatother \loadmsbm

\def\R{\Bbb R}

\def\Z{\Bbb Z}

\def\no{\noindent}
\def\f{\frac}

\def\pa{\partial}

\def\na{\nabla}

\def\al{\alpha}

\def\endproof{\hphantom{MM}\hfill\llap{$\square$}\goodbreak}

\newcommand{\beq}{\begin{equation}}
\newcommand{\eeq}{\end{equation}}
\newcommand{\ben}{\begin{eqnarray}}
\newcommand{\een}{\end{eqnarray}}
\newcommand{\beno}{\begin{eqnarray*}}
\newcommand{\eeno}{\end{eqnarray*}}

\newtheorem{Theorem}{Theorem}[section]
\newtheorem{Lemma}[Theorem]{Lemma}
\newtheorem{Def}{Definition}[section]
\newtheorem{Remark}{Remark}[section]

\begin{document}
\title{On the regularity criterion of weak solution for the 3D viscous Magneto-hydrodynamics equations}

\author{Qionglei Chen $^\dag$,
Changxing Miao $^\dag$,
Zhifei Zhang $^\ddag$\\
{\small  $^\dag$Institute of Applied Physics and Computational Mathematics,}\\
    {\small P.O. Box 8009, Beijing 100088, P. R. China.}\\
     {\small (chen\_qionglei@iapcm.ac.cn  and   miao\_changxing@iapcm.ac.cn)}\\
{\small $^\ddag$ School of Mathematical Science, Peking University,}\\
{\small Beijing 100871, P. R.  China.}\\
{\small(zfzhang@math.pku.edu.cn)  }}
\date{}
\maketitle

\vspace{-1.2in} \vspace{.9in} \vspace{0.2cm}
{\bf Abstract.} We improve and extend some known
regularity criterion of weak solution for the 3D viscous Magneto-hydrodynamics equations
by means of the Fourier localization technique and Bony's para-product decomposition.

\vspace{0.2cm}
{\bf Key words.} MHD equations, Weak solution, Regularity criterion,  Fourier Localization, Bony's para-product
decompoition.
 \vspace{0.2cm}

{\bf AMS subject classifications.} 76W05 35B65

\vskip0.2cm

\section{Introduction}
\vskip0.2cm

In this paper, we consider the 3D incompressible magneto-hydrodynamics(MHD) equations
\begin{align} \label{1.1}(\rm MHD)\,\, \left\{
\begin{aligned}
&\frac{\pa u}{\pa t}-\nu\Delta u+u\cdot \nabla u=-\nabla p-\frac{1}{2}\nabla b^2+b\cdot \nabla
b,\\
&\frac{\pa b}{\pa t}-\eta\Delta b+u\cdot \nabla b=b\cdot \nabla u ,\\
&\nabla\cdot u=\nabla\cdot b=0,\\
&u(0,x)=u_0(x),\quad b(0,x)=b_0(x).
\end{aligned}
\right. \end{align}
Here $u$, $b$ describe
the flow velocity vector and the magnetic field vector
respectively, $p$ is a scalar pressure, $\nu>0$ is the kinematic viscosity, $\eta>0$ is the magnetic diffusivity, while $u_0$ and $b_0$ are
 the given initial velocity and initial magnetic field
with $\nabla\cdot u_0=\nabla\cdot b_0=0$. If $\nu=\eta=0$, (\ref{1.1}) is called the ideal MHD equations.

As same as the 3D Navier-Stokes equations, the regularity of weak solution for the 3D MHD equations remains open\cite{ST}.
For the 3D Navier-Stokes equations, the Serrin-type criterion states that a Leray-Hopf weak
solution $u$ is regular provided the following condition holds\cite{Bei,ES, Gi,Ser,Wa}:
\ben
&u\in L^q(0,T; L^p(\R^3)), \quad\mbox{for}\quad \frac 2 q+\frac 3 p\le 1,\quad 3\le p\le\infty,\label{1.2}\een
or

\ben\na u\in L^q(0,T; L^p(\R^3)), \quad\mbox{for}\quad
\frac 2 q+\frac 3 p\le 2,\quad \f 32<p\le\infty.\label{1.3}
\een
Recently, Chen and Zhang \cite{CZ} have refined the above conditions
as follows: If there exists a small $\varepsilon_0$ such that for
any $t\in (0,T), u$ satisfies
\beq \lim_{\varepsilon\rightarrow
0}\sup_j2^{js
q}\int_{t-\varepsilon}^t\|\Delta_ju(\tau)\|_p^qd\tau\le
\varepsilon_0,\label{1.4}
\eeq
with $\f 2 q+\f 3 p =1+s, \f 3
{1+s}<p\le \infty, -1<s\le 1$, and $(p,s)\neq (\infty,1)$, then $u$
is regular in $(0,T]\times \R^3$, where $\Delta_j$ denotes the frequency localization operator.
 For the marginal case($p=3 ,q=\infty$),
Cheskidov and Shvydkoy \cite{Che} have refined \eqref{1.2} to
\beq
u\in C([0,T]; B^{-1}_{\infty,\infty}).\label{1.5}
\eeq Here $B^{-1}_{\infty, \infty}$ stands for the inhomogenous Besov spaces,
see Section 2 for the definitions.

Wu \cite{Wu1,Wu2} extended some Serrin-type criteria for the
Navier-Stokes equations to the MHD equations imposing conditions
on both the velocity field $u$ and the magnetic field $b$.
However, some numerical experiments \cite{PP} seem to indicate
that the velocity field
 plays the more important role than the magnetic field
in the regularity theory of solutions to the MHD equations.
Recently, He, Xin\cite{HX}, and Zhou\cite{Zhou1} have proved some
regularity criteria to the MHD equations which do not impose any
condition on the magnetic field $b$. Precisely, they showed that
the weak solution remains smooth on $(0,T]\times \R^3$ if the
velocity $u$ satisfies  one of the following conditions
\begin{align}\label{1.6}
&u\in L^q(0,T; L^p(\R^3)),\quad \frac{2}{q}+\frac{3}{p}\le1,\quad
3<p\le\infty;
\\\label{1.7}
& u\in C([0,T]; L^3(\R^3));\\\label{1.8}
&\na u\in L^q(0,T; L^p(\R^3)),\quad \frac{2}{q}+\frac{3}{p}\le2,
\quad \frac32<p\le\infty.
\end{align}
Meanwhile, inspired by the pioneering work of Constantin
and Fefferman  \cite{Con} where the regularity condition of the
direction of vorticity was used to describe the regularity
criterion to the Navier-Stokes equations, He and Xin \cite{HX}
showed that the weak solution remains smooth on $(0,T]\times \R^3$
 if the vorticity of the velocity
$w=\na\times u$ satisfies  the following condition
\begin{align}\label{1.9}
\big|w(x+y, t)-w(x, t)\big|\le K|w(x+y,t)||y|^{\frac12}\quad \mbox{if} \,\,\,|y|\le\rho\quad|w(x+y,t)|\ge\Omega,
\end{align}
for $t\in[0,T]$ and three positive constants $K$, $\rho$,
$\Omega$.

For the marginal case $p=\infty$ in \eqref{1.8},
Chen, Miao and Zhang \cite{CMZ} proved a Beale-Kato-Majda  criterion in terms of the vorticity of the velocity $u$ only
by means of  the  Littlewood-Paley
decomposition.

For the generalized MHD equations with fractional dissipative effect, Wu\cite{Wu3, Wu4}  established some regularity results
in terms of the velocity only.

The purpose of this paper is to improve and extend some known regularity criterion
of weak solution for the MHD equations by  means of the Fourier localization technique and
Bony's para-product decomposition \cite{Bon, Ch}.
Let us firstly recall the definition of weak solution.
\begin{Def}
The vector-valued function $(u,b)$ is called a weak solution of (\ref{1.1}) on $(0,T)\times \R^3$ if
it satisfies the following conditions:
\begin{flushleft}
\textrm{(1)}\, $(u,b)\in L^\infty(0,T;L^2(\R^3))\cap L^2(0,T;H^1(\R^3))$;

\textrm{(2)}\, $\textrm{div}\, u=\textrm{div}\, b=0$ in the sense of distribution;

\textrm{(3)}\, For any function $\psi(t,x)\in C^\infty_0((0,T)\times \R^3)$ with $\textrm{div} \psi=0$, there hold
$$
\int_0^T\int_{\R^3}\{u\cdot \psi_t-\nu\na u\cdot \na \psi+\na \psi:(u\otimes u-b\otimes b)\}dxdt=0,
$$
and
$$
\int_0^T\int_{\R^3}\{b\cdot \psi_t-\eta\na b\cdot \na \psi+\na
\psi:(u\otimes b-b\otimes u)\}dxdt=0.
$$
\end{flushleft}
\end{Def}
Similar to the Navier-Stokes equation, the global existence of weak
solutions to the MHD equations can be proved by using the Galerkin's
method and compact argument, see\cite{Duv}. Now we state our main
result as follows.

\begin{Theorem}\label{main}
Let $(u_0,b_0)\in L^2(\R^3)$  with $\na\cdot u_0=\na\cdot
b_0=0$. Assume that $(u,b)$ is a weak solution to (\ref{1.1}) on $(0,T)\times \R^3$ with $0\le T\le \infty$.
If  the velocity $u(t)$  satisfies
\begin{align}\label{1.10}
u(t)\in L^q(0,T; B^s_{p,\infty}),
\end{align}
with $\f 2 q+\f 3 p =1+s, \f 3 {1+s}<p\le \infty, -1<s\le 1$, and
$(p,s)\neq (\infty,1)$. Then the solution $(u,b)$ is regular on $(0,T]\times \R^3$.
\end{Theorem}

\begin{Remark}
By the embedding $L^p \subsetneq B^0_{p,\infty}$, we see that our
result is an improvement of (\ref{1.6}) and (\ref{1.8}). In
addition, we establish the regularity criterion of weak solution for the MHD equation in the framework of
Besov spaces with negative index in terms of the velocity only.
On the other hand, the method in this paper can be applied  to the generalized MHD equations, please refer to \cite{Wu3, Wu4}
for details.\end{Remark}

\begin{Remark}
In the case of $s=0$ or $s=1$, Kozono,  Ogawa  and Taniuchi\cite{Koz}
proved the similar results for the Navier-Stokes equations by using the Logarithmic Sobolev
inequality in the Besov spaces. However, if we try to use their
method to our case, we can only obtain the regularity criterion in
terms of both the velocity field  $u$ and the magnetic field $b$.
\end{Remark}

\begin{Remark}
Chen, Miao and Zhang\cite{CMZ} proved the marginal case $(p,s)=
(\infty,1)$ by using a different argument. However, the method of
\cite{CMZ} can't also be applied to the present case.
\end{Remark}

\begin{Remark}
The regularity of weak solution $(u,b)$ under the condition
\begin{equation}\label{1.11}
u\in C(0,T; B^{-1}_{\infty,\infty})
\end{equation}
remains unknown.  One easily checks that it is  the special case of
the endpoint case of \eqref{1.10} in Theorem \ref{main} with  $s=-1$.

\end{Remark}

\noindent{\bf Notation:} Throughout the paper, $C$ stands for   a generic  constant.
We  will use the notation $A\lesssim B$ to denote the relation  $A\le CB$ and
the notation $A\approx B$ to denote the relations  $A\lesssim B$ and $B\lesssim A$.
Further, $\|\cdot\|_{p}$ denotes the norm of the Lebesgue space $L^p$.

\vskip0.2cm

\section{Preliminaries}

\vskip0.2cm

In this section, we are going to recall some basic facts on
Littlewood-Paley  theory, one may check \cite{Ch} for more details.

Let ${\cal S}(\R^3)$ be
the Schwartz class of rapidly decreasing functions. Given $f\in
{\cal S}(\R^3)$, its Fourier transform ${\cal F}f=\hat f$ is
defined by
$$
\hat f(\xi)=(2\pi)^{-\frac{3}2}\int_{\R^3}e^{-ix\cdot \xi}f(x)dx.
$$
Choose two nonnegative radial functions $\chi$, $\varphi \in {\cal
S}(\R^3)$ supported respectively in ${\cal B}=\{\xi\in\R^3,\,
|\xi|\le\frac{4}{3}\}$ and ${\cal C}=\{\xi\in\R^3,\,
\frac{3}{4}\le|\xi|\le\frac{8}{3}\}$ such that
\beno
\chi(\xi)+\sum_{j\ge0}\varphi(2^{-j}\xi)=1,\quad\xi\in\R^3.\\
\eeno
Set $\varphi_j(\xi)=\varphi(2^{-j}\xi)$ and
let $h={\cal F}^{-1}\varphi$ and $\tilde{h}={\cal F}^{-1}\chi$.
Define the frequency localization operators:
\begin{align}
&\Delta_jf=\varphi(2^{-j}D)f=2^{3j}\int_{\R^3}h(2^jy)f(x-y)dy,\quad\mbox{for}\quad
j\geq 0, \nonumber\\
&S_jf=\chi(2^{-j}D)f=\sum_{-1\leq k\le
j-1}\Delta_kf=2^{3j}\int_{\R^3}\tilde{h}(2^jy)f(x-y)dy,\quad\mbox{and}\nonumber\\
&\Delta_{-1}f=S_{0}f, \qquad \qquad\Delta_{j}f=0\quad\mbox{for}\quad
j\le -2. \label{2.1}\end{align} Formally, $\Delta_j=S_{ j+1}-S_{j}$
is a frequency projection into the annulus $\{|\xi|\approx 2^j\}$,
and  $S_j$ is a frequency projection into the ball $\{|\xi|\lesssim
2^j\}$. One easily verifies that with the above choice of $\varphi$
\begin{align}\label{2.2}
\Delta_j\Delta_kf\equiv0\quad i\!f\quad|j-k|\ge 2\quad and
\quad \Delta_j(S_{k-1}f\Delta_k
f)\equiv0\quad i\!f\quad|j-k|\ge 5.
\end{align}
We now introduce the following definition of inhomogenous Besov spaces by means of  Littlewood-Paley projection  $\Delta_j$ and $S_j$:
\begin{Def}
Let $s\in \R, 1\le p,q\le\infty$, the inhomogenous Besov space
$B^s_{p,q}$ is defined by
$$B^s_{p,q}=\{f\in {\cal S}'(\R^3); \quad \|f\|_{B^s_{p,q}}<\infty\},$$ where
$$\|f\|_{B^s_{p,q}}=\left\{\begin{array}{l}
\displaystyle\bigg(\sum_{j= -1}^{\infty}2^{jsq}\|\Delta_j
f\|_{L^p}^q\bigg)^{\frac 1
q},\quad \hbox{for}\quad q<\infty,\\
\displaystyle\sup_{j\geq -1}2^{js}\|\Delta_jf\|_{L^p}, \quad
\hbox{ for} \quad q=\infty.
\end{array}\right.
$$
\end{Def}
Let us point out that $B^s_{2,2}$ is the usual Sobolev space
$H^s$ and that $B^{s}_{\infty,\infty}$ is the usual H\"{o}lder
space $C^{s}$ for $s\in \R\setminus\Z$. We refer to \cite{Tri} for more details.

We now recall the para-differential calculus which enables us to define a generalized
product between distributions, which is continuous in many functional spaces where the usual
product does not make sense (see \cite{Bon}). The paraproduct between $u$ and $v$ is defined
by
\begin{align}\label{2.3}
T_uv\triangleq\sum_{j}S_{j-1}u\Delta_jv.
\end{align}
Formally, we have the following Bony's decomposition:
\beq\label{2.4}
uv=T_uv+T_vu+R(u,v),
\eeq
with $$R(u,v)=\sum_{|j'-j|\le 1}\Delta_ju \Delta_{j'}v,$$
and we also denote $$T'_{u}v\triangleq T_{u}v+R(u,v).$$

Let us conclude this section by recalling the Bernstein's inequality
which will be frequently used in the proof of Theorem 1.1.

\begin{Lemma}\cite{Ch}\label{Lem2.1}
Let $1\le p\le q\le\infty$. Assume that $f\in L^p$, then there
exists a constant $C$ independent of $f$, $j$ such that
 \ben
{\rm supp}\hat{f}\subset\{|\xi|\leq C2^j\}\Longrightarrow
\|\partial^\alpha f\|_{L^q}\le
C3^{j{|\alpha|}+3j(\frac{1}{p}-\frac{1}{q})}\|f\|_{L^p},\label{2.5}
\een   \ben{\rm supp}\hat{f}\subset\Big\{\f1{C}2^j\leq |\xi|\leq C
2^j\Big\}\Longrightarrow \|f\|_{L^p}\le
C2^{-j|\alpha|}\sup_{|\beta|=|\al|}\|\partial^\beta
f\|_{L^p}.\label{2.6} \een
\end{Lemma}
\vskip0.2cm

\section{Proof of Theorem 1.1}
\vskip0.2cm

Since the weak solution $(u(t),b(t))\in L^2(0,T;H^1(\R^3))$, for any time
interval $(0,\delta)$, there exists an $\varepsilon\in (0,\delta)$
such that
$(u(\varepsilon), b(\varepsilon))\in H^1(\R^3)$. It is well
known that there exist a maximal existence time $T_0>0$ and a
unique strong solution
$$(\widetilde{u}(t),\widetilde{b}(t))\in {\cal
X}(\varepsilon,  T_0)\triangleq C([\varepsilon,T_0);H^1(\R^3))\cap
C^1((\varepsilon,T_0);H^{1}(\R^3))\cap C((\varepsilon,T_0);
H^{3}(\R^3))$$ which is the same as the weak solution $(u,b)$ on
$(\varepsilon, T_0)\cite{Duv,ST}$. In order to complete the proof of
Theorem 1.1, it suffices to show that the strong solution
$(u(t),b(t))$ can be extended after $t=T_0$ in the class ${\cal
X}(\varepsilon,  T_0)$ under the condition of Theorem 1.1. For the
convenience, we set $\nu=\eta=1$ and $\varepsilon=0$ in what
follows. We denote
$$u_k=\Delta_ku,\quad b_k=\Delta_kb,\quad
\pi_k=\Delta_k\pi,$$
here $\pi=p+\frac12 b^2$. We get by applying the operation $\Delta_k$ to
both sides of \eqref{1.1} that
\begin{align} \label{3.1}
\left\{
\begin{aligned}
&{\pa_t u_k}-\Delta u_k+\Delta_k(u\cdot \nabla u)-\Delta_k(b\cdot \nabla
b)=-\nabla\pi_k,\\
&{\pa_t b_k}-\Delta b_k+\Delta_k(u\cdot \nabla b)-\Delta_k(b\cdot \nabla u)=0.
\end{aligned}
\right.
\end{align}
Multiplying the first equation of (\ref{3.1}) by $u_k$ and  the second one of (\ref{3.1}) by $b_k$,
we obtain by Lemma 2.1 for $k\ge 0$ that
\begin{align}
&\frac{1}{2}\frac{d}{dt}\|u_k(t)\|_2^2+c2^{2k}\|
u_k(t)\|_2^2=-\big<\Delta_k(u\cdot\na u),\,u_k\big>+\big<\Delta_k(b\cdot\na b),\,u_k\big>,\label{3.2}\\
&\frac{1}{2}\frac{d}{dt}\|b_k(t)\|_2^2+c2^{2k}\| b_k(t)\|_2^2=
-\big<\Delta_k(u\cdot\na b),\,b_k\big>+\big<\Delta_k(b\cdot\na u),\,b_k\big>.\label{3.3}
\end{align}
Set
$$F_k(t)\triangleq\bigl(\|u_k(t)\|_2^2+\|b_k(t)\|_2^2\bigr)^\f12.$$
Then we get by adding (\ref{3.2}) and (\ref{3.3}) that
\begin{align}
\frac{1}{2}\frac{d}{dt}F_k(t)^2+c2^{2k}F_k(t)^2
=&-\big<\Delta_k(u\cdot\na u),\,u_k\big>+\big<\Delta_k(b\cdot\na b),\,u_k\big>\nonumber\\
&-\big<\Delta_k(u\cdot\na b),\,b_k\big>+\big<\Delta_k(b\cdot\na u),\,b_k\big>.\label{3.4}
\end{align}
Noting that
$$\big<u\cdot\na u_k, u_k\big>=\big<u\cdot\na b_k, b_k\big>=0,$$
$$\big<b\cdot\na b_k, u_k\big>+\big<b\cdot\na u_k, b_k\big>=0.$$
The right hand side of \eqref{3.4} can be written as
\begin{align}
&\big<[u, \Delta_k]\cdot\na u, u_k\big>-\big<[b, \Delta_k]\cdot\na b, u_k\big>+\big<[u, \Delta_k]\cdot\na b, b_k\big>
-\big<[b, \Delta_k]\cdot\na u, b_k\big>\nonumber\\ &\triangleq I+II+III+IV, \nonumber
\end{align}
where $[A,B]\triangleq AB-BA$. Using the Bony's decomposition (\ref{2.4}), we  rewrite $I$ as
\begin{align*}
I&=\big<[T_{u^i}, \Delta_k]\pa_i u, u_k\big>
+\big<T'_{\Delta_k\pa_i u}u^i, u_k\big>-\big<\Delta_kT_{\pa_i u}u^i, u_k\big>-\big<\Delta_kR(u^i,\pa_i u), u_k\big>\\
&\triangleq I_1+I_2+I_3+I_4.
\end{align*}
In view of the support of the Fourier transform of the term
$T_{\pa_i u}u^i$, we have
$$
\Delta_k T_{\pa_i u}u^i=\sum_{|k'-k|\le4}\Delta_k(S_{k'-1}(\pa_i
u)u^i_{k'}).
$$
This helps us to get by Lemma \ref{Lem2.1}
\beq
|I_3| \lesssim \|u_k\|_2\sum_{|k'-k|\le4}\|\na
S_{k'-1}u\|_\infty\|u_{k'}\|_2.\label{3.5}
\eeq
Since $\textrm{div}\,u=0$, we have
$$
\Delta_kR(u^i,\pa_i u)=\sum_{k',k''\ge k-2;|k'-k''|\le 1}\pa_i\Delta_k(\Delta_{k'}u^i\Delta_{k''}u).
$$
This together with Lemma 2.1 yields
\begin{align}
|I_4| \lesssim 2^k\|u_k\|_\infty\sum_{k'\ge k-2}\|u_{k'}\|_2^2.\label{3.6}
\end{align}
Using  the definition of $T'_{\Delta_k\pa_i u}u^i$, we have
\begin{align}
T'_{\Delta_k\pa_i u}u^i=\sum_{k'\ge k-2} S_{k'+2}\Delta_k\pa_i
u\Delta_{k'}u^i.\nonumber
\end{align}
Note that $S_{k'+2}\Delta_k u=\Delta_ku$ for $k'>k$, we get
$$
I_2=\sum_{k-2\le k'\le k}\big\langle S_{k'+2}\Delta_k\pa_i
u\Delta_{k'}u^i,u_k\big\rangle,
$$
from which and Lemma 2.1, it follows that
\begin{align}\label{3.7}
I_2&\lesssim \|u_k\|_2\sum_{|k'-k|\le2}\|\na
S_{k'-1}u\|_\infty\|u_{k'}\|_2+2^k\|u_k\|_\infty\sum_{k'\ge
k-2}\|u_{k'}\|_2^2.
\end{align}
Making use of  the definition of $\Delta_k$, we have
\begin{align*}\nonumber
[T_{u^i}, \Delta_k]\pa_i u&=\sum_{|k'-k|\le4}[S_{k'-1}u^i, \Delta_k]\pa_i u_{k'}\\
&=\sum_{|k'-k|\le4}2^{3k}\int_{\R^3}h(2^k(x-y))(S_{k'-1}u^i(x)-S_{k'-1}u^i(y))\pa_i u_{k'}(y)dy\\
&=\sum_{|k'-k|\le4}2^{4k}\int_{\R^3}\int_0^1y\cdot\na S_{k'-1}u^i(x-\tau y)d\tau\pa_i h(2^ky)u_{k'}(x-y)dy,
\end{align*}
from which and the Minkowski inequality, we infer that
\begin{align}\label{3.8}
|I_1| \lesssim \|u_k\|_2\sum_{|k'-k|\le4}\|\na
S_{k'-1}u\|_\infty\|u_{k'}\|_2.
\end{align}
By summing up (\ref{3.5})-(\ref{3.8}), we obtain
\begin{align}
|I| \lesssim \|u_k\|_2\sum_{|k'-k|\le4}\|\na
S_{k'-1}u\|_\infty\|u_{k'}\|_2+2^k\|u_k\|_\infty\sum_{k'\ge
k-2}\|u_{k'}\|_2^2.\label{3.9}
\end{align}

Similar arguments as in deriving \eqref{3.9} can be used to get that
\begin{align}
|II+IV| \lesssim &\|u_k\|_\infty\sum_{|k'-k|\le4}\|\na
S_{k'-1}b\|_2\|b_{k'}\|_2+2^k\|u_k\|_\infty\sum_{k'\ge
k-2}\|b_{k'}\|_2^2\nonumber\\
&+\|b_k\|_2\sum_{|k'-k|\le4}(\|\na
S_{k'-1}u\|_\infty\|b_{k'}\|_2+\|\na
S_{k'-1}b\|_2\|u_{k'}\|_\infty)\nonumber\\&+2^k\|b_k\|_2\sum_{k',k''\ge
k-2;|k'-k''|\le 1}\|u_{k'}\|_\infty\|b_{k''}\|_2, \label{3.10}
\end{align}
and
\begin{align}
|III| \lesssim &\|b_k\|_2\sum_{|k'-k|\le4}(\|\na
S_{k'-1}u\|_\infty\|b_{k'}\|_2+\|\na
S_{k'-1}b\|_2\|u_{k'}\|_\infty)\nonumber\\&+2^k\|b_k\|_2\sum_{k',k''\ge
k-2;|k'-k''|\le 1}\|u_{k'}\|_\infty\|b_{k''}\|_2.\label{3.11}
\end{align}

Combining \eqref{3.9}-\eqref{3.11} with \eqref{3.4},  we easily get
for $k\ge 0$ that
\begin{align}\label{3.13}
&\f d {dt}F_k(t)^2+2^{2k}F_k(t)^2\nonumber\\
&\lesssim (\|u_k\|_2+\|b_k\|_2)\sum_{|k'-k|\le4}\|\na S_{k'-1}u\|_\infty(\|u_{k'}\|_2+\|b_{k'}\|_2)\nonumber\\
&\quad+\sum_{|k'-k|\le4}\|\na S_{k'-1}b\|_2(\|b_{k'}\|_2\|u_k\|_\infty+\|u_{k'}\|_\infty\|b_k\|_2)\nonumber\\
&\quad+2^k\|b_k\|_2\sum_{k',k''\ge k-2;|k'-k''|\le
1}\|u_{k'}\|_\infty\|b_{k''}\|_2\nonumber\\
&\quad+2^k\|u_k\|_\infty\sum_{k'\ge
k-2}(\|b_{k'}\|_2^2+\|u_{k'}\|_2^2).
\end{align}
Making use of $ B^{\frac2q+\frac3p-1}_{p,\infty}(\R^3)\hookrightarrow B^{\frac2q-1}_{\infty,\infty}(\R^3)$,
we only need to deal with the case when $p=+\infty$ since the other cases can be deduced from it by above Sobolev embedding.
Here we omit the details. By the restrictions on $p,q,s$, we see that $s=\f2q-1$ and $q\in (1,+\infty)$.
\vspace{0.3cm}

\no\textbf{Case 1}. $q \in (1,2]$.

\vspace{0.2cm}

Integrating (\ref{3.13}) with respect to $t$, we deduce that
\begin{align}
&F_k(t)^2-F_k(0)^2+2^{2k}\int_0^tF_k(\tau)^2d\tau\nonumber\\
&\lesssim \int_0^t(\|u_k\|_2+\|b_k\|_2)\sum_{|k'-k|\le4}\|\na S_{k'-1}u\|_\infty(\|u_{k'}\|_2+\|b_{k'}\|_2)d\tau\nonumber\\
&\quad+\int_0^t\sum_{|k'-k|\le4}\|\na S_{k'-1}b\|_2(\|b_{k'}\|_2\|u_k\|_\infty+\|u_{k'}\|_\infty\|b_k\|_2)d\tau\nonumber\\
&\quad+\int_0^t2^k\|b_k\|_2\sum_{k',k''\ge k-2;|k'-k''|\le
1}\|u_{k'}\|_\infty\|b_{k''}\|_2 d\tau\nonumber\\
&\quad+\int_0^t2^k\|u_k\|_\infty\sum_{k'\ge k-2}(\|b_{k'}\|_2^2+\|u_{k'}\|_2^2)d\tau\nonumber\\
&\triangleq \Pi_1+\Pi_2+\Pi_3+\Pi_4.\label{3.14}
\end{align}
Take $\rho \in (\f12, 1)$ and set
\begin{align*}
A(t)\triangleq\sup_{k\ge-1} 2^{k\rho}F_k(\tau),\quad
B(t)=\sup_{k\ge-1}2^{2k(\rho+1)}\int_0^tF_k(\tau)^2d\tau.
\end{align*}
We get by using Lemma 2.1 that
\begin{align}
2^{2k\rho}\Pi_1 &\lesssim
\int_0^t2^{2k\rho}(\|u_{k}\|_2+\|b_{k}\|_2)\sum_{|k'-k|\le4}\sum_{k''=-1}^{k'-2}\|u_{k''}\|_\infty2^{k''}
(\|u_{k'}\|_2+\|b_{k'}\|_2)d\tau\nonumber\\
&\lesssim
\int_0^t\|u(\tau)\|_{B^s_{\infty,\infty}}A(\tau)2^{2k\rho}(\|u_{k}\|_2+\|b_{k}\|_2)
\sum_{|k'-k|\le4}2^{-k'\rho}\sum_{k''=-1}^{k'-2}2^{k''(2-\frac 2 q)}d\tau\nonumber\\
&\lesssim
\int_0^t\|u(\tau)\|_{B^s_{\infty,\infty}}A(\tau)2^{k(\rho +2-\frac 2 q)}F_k(\tau) d\tau\nonumber\\
&\lesssim  \Bigl(\int_0^t\|u(\tau)\|_{B^s_{\infty,\infty}}^qA(\tau)^2d\tau\Bigr)^\f1qB(t)^{1-\f1q},\label{3.15}
\end{align}
where we used the fact that $1<q\le 2$ in the last two inequalities.
Similarly, we get by using Lemma 2.1 and the fact that $\rho<1$ and
$1<q\le 2$ that
\begin{align}
2^{2k\rho}\Pi_2 &\lesssim
\int_0^t\|u(\tau)\|_{B^s_{\infty,\infty}}2^{k(1-\frac2q+2\rho)}\sum_{|k'-k|\le4}\sum_{k''=-1}^{k'-2}\|b_{k''}\|_22^{k''}\big(\|b_{k'}\|_2
+\|b_{k}\|_2
\big)d\tau\nonumber\\
&\lesssim
\int_0^t\|u(\tau)\|_{B^s_{\infty,\infty}}A(\tau)2^{k(1-\frac2q+2\rho)}
\sum_{|k'-k|\le4}2^{k'(1-\rho)}\big(\|b_{k'}\|_2+\|b_{k}\|_2\big)d\tau\nonumber\\
&\lesssim \Bigl(\int_0^t\|u(\tau)\|_{B^s_{\infty,\infty}}^qA(\tau)^2d\tau\Bigr)^\f1qB(t)^{1-\f1q},\label{3.16}\\
2^{2k\rho}\Pi_3
&\lesssim \int_0^t\|u(\tau)\|_{B^s_{\infty,\infty}}2^{k(2\rho+1)}\|b_k\|_2\sum_{k'\ge k-2}\|b_{k'}\|_22^{k'(1-\frac2q)}d\tau\nonumber\\
&\lesssim \int_0^t\|u(\tau)\|_{B^s_{\infty,\infty}}A(\tau)2^{k(2\rho+1)}\|b_k\|_2\sum_{k'\ge k-2}2^{k'(1-\frac2q-\rho)}d\tau\nonumber\\
&\lesssim
\Bigl(\int_0^t\|u(\tau)\|_{B^s_{\infty,\infty}}^qA(\tau)^2d\tau\Bigr)^\f1qB(t)^{1-\f1q},\label{3.17}\\
2^{2k\rho}\Pi_4 &\lesssim
\int_0^t\|u(\tau)\|_{B^s_{\infty,\infty}}2^{k(2\rho+2-\f2q)}\sum_{k'\ge
k-2}(\|b_{k'}\|_2^2+\|u_{k'}\|_2^2)
d\tau\nonumber\\
&\lesssim \Bigl(\int_0^t\|u(\tau)\|_{B^s_{\infty,\infty}}^qA(\tau)^2d\tau\Bigr)^\f1qB(t)^{1-\f1q}.\label{3.18}
\end{align}

On the other hand, the strong solution $(u,b)$ also satisfies the energy equality
$$
\|u(t)\|_2^2+\|b(t)\|_2^2+2\int_0^t(\|\na u(\tau)\|_2^2+\|\na b(\tau)\|^2_2)d\tau=\|u_0\|_2^2+\|b_0\|_2^2,
$$
hence, we have
\ben
\|\Delta_{-1}u(t)\|_2+\|\Delta_{-1}b(t)\|_2\le C(\|u_0\|_2+\|b_0\|_2).\label{3.19}
\een
Thus, summing up \eqref{3.14}-\eqref{3.19}, we get by the Young's inequality that
\beno
A(t)^2 \le C(A(0)^2+\|u_0\|_2^2+\|b_0\|_2^2)+C\int_0^t\|u(\tau)\|_{B^s_{\infty,\infty}}^qA(\tau)^2d\tau,
\eeno
which together with the Gronwall inequality yields that
\ben
A(t)^2\le C(A(0)^2+\|u_0\|_2^2+\|b_0\|_2^2)
\exp\Bigl(C\int_0^t\|u(\tau)\|_{B^s_{\infty,\infty}}^qd\tau\Bigr).\label{3.20}
\een

\vspace{0.3cm}
\no\textbf{Case 2}. $q\in (2,+\infty)$
\vspace{0.2cm}

Multiplying (\ref{3.13}) by $2^{2k(q-1)}F_k(t)^{2(q-1)}$, then integrating the resulting equation with respect to $t$
leads to the result  for $k\ge 0$
\begin{align}\label{3.21}
&2^{2k(q-1)}F_k(t)^{2q}-2^{2k(q-1)}F_k(0)^{2q}+2^{2kq}\int_0^tF_k(\tau)^{2q}d\tau\nonumber\\
&\lesssim \int_0^t2^{2k(q-1)}F_k(\tau)^{2(q-1)}(\|u_k\|_2+\|b_k\|_2)\sum_{|k'-k|\le4}\|\na S_{k'-1}u\|_\infty(\|u_{k'}\|_2+\|b_{k'}\|_2)d\tau\nonumber\\
&\quad+\int_0^t2^{2k(q-1)}F_k(\tau)^{2(q-1)}\sum_{|k'-k|\le4}\|\na S_{k'-1}b\|_2(\|b_{k'}\|_2\|u_k\|_\infty+\|u_{k'}\|_\infty\|b_k\|_2)d\tau\nonumber\\
&\quad+\int_0^t2^{2k(q-1)}F_k(\tau)^{2(q-1)}2^k\|b_k\|_2\sum_{k',k''\ge
k-2;|k'-k''|\le1}\|u_{k'}\|_\infty\|b_{k''}\|_2 d\tau\nonumber\\
&\quad+\int_0^t2^{2k(q-1)}F_k(\tau)^{2(q-1)}2^k\|u_k\|_\infty\sum_{k'\ge k-2}(\|b_{k'}\|_2^2+\|u_{k'}\|_2^2)d\tau\nonumber\\
&\triangleq \Pi_1+\Pi_2+\Pi_3+\Pi_4.
\end{align}
Set
\begin{align*}
&A(t)\triangleq\sup_{k\ge-1} 2^{(1-\frac{1}{q})k}F_k(t),\qquad
 B(t)\triangleq\sup_{k\ge -1}2^{2kq}\int_0^tF_k(\tau)^{2q}d\tau.
\end{align*}
We obtain that by  Lemma 2.1
\begin{align}
\Pi_1&\lesssim \int_0^t2^{2k(q-1)}F_k(\tau)^{2(q-1)}(\|u_{k}\|_2+\|b_{k}\|_2)\sum_{|k'-k|\le4}\sum_{k''=-1}^{k'-2}\|u_{k''}\|_\infty2^{k''}
(\|u_{k'}\|_2+\|b_{k'}\|_2)d\tau\nonumber\\
&\lesssim \int_0^t2^{2k(q-1)}F_k(\tau)^{2(q-1)}\|u(\tau)\|_{B^s_{\infty,\infty}}A(\tau)^2d\tau\nonumber\\
&\lesssim  \Bigl(\int_0^t\|u(\tau)\|_{B^s_{\infty,\infty}}^qA(\tau)^{2q}d\tau\Bigr)^\f1qB(t)^\f{q-1}q.\label{3.22}
\end{align}
Similarly, we get by using Lemma 2.1 and the fact that $q<+\infty$ that
\begin{align}
\Pi_2
&\lesssim \int_0^t2^{2k(q-1)}F_k(\tau)^{2(q-1)}\|u(\tau)\|_{B^s_{\infty,\infty}}A(\tau)2^{k(1-\frac2q)}\sum_{|k'-k|\le4}\sum_{k''=-1}^{k'-2}2^{\frac{k''}q}\big(\|b_{k'}\|_2
+\|b_{k}\|_2\big)d\tau\nonumber\\
&\lesssim \int_0^t2^{2k(q-1)}F_k(\tau)^{2(q-1)}\|u(\tau)\|_{B^s_{\infty,\infty}}A(\tau)\sum_{|k'-k|\le4}2^{k(1-\frac1q)}\big(\|b_{k'}\|_2
+\|b_{k}\|_2\big)d\tau\nonumber\\
&\lesssim \Bigl(\int_0^t\|u(\tau)\|_{B^s_{\infty,\infty}}^qA(\tau)^{2q}d\tau\Bigr)^\f1qB(t)^\f{q-1}q,\label{3.23}\end{align}

\begin{align}\nonumber
\Pi_3\lesssim& \int_0^t2^{2k(q-1)}F_k(\tau)^{2(q-1)}\|u(\tau)\|_{B^s_{\infty,\infty}}\sum_{k'\ge k-2}\|b_{k'}\|_22^{k'(1-\frac2q)}2^{k}\|b_k\|_2d\tau\nonumber\\
\lesssim&   \int_0^t2^{2k(q-1)}F_k(\tau)^{2(q-1)}\|u(\tau)\|_{B^s_{\infty,\infty}}A(\tau)\sum_{k'\ge k-2}2^{-\frac{k'}q}2^{k}\|b_k\|_2d\tau\nonumber\\
\lesssim&
\Bigl(\int_0^t\|u(\tau)\|_{B^s_{\infty,\infty}}^qA(\tau)^{2q}d\tau\Bigr)^\f1qB(t)^\f{q-1}q,\label{3.24}\\
\Pi_4
\lesssim &\int_0^t2^{2k(q-1)}F_k(\tau)^{2(q-1)}\|u(\tau)\|_{B^s_{\infty,\infty}}\nonumber\\
&\quad\times \sum_{k'\ge k-2}(\|b_{k'}\|_2^2+\|u_{k'}\|_2^2)2^{k'(2-\frac2q)}2^{(k-k')(2-\frac2q)}d\tau
\nonumber\\
\lesssim &\int_0^t2^{2k(q-1)}F_k(\tau)^{2(q-1)}\|u(\tau)\|_{B^s_{\infty,\infty}}A(\tau)^2d\tau\nonumber\\
\lesssim & \Bigl(\int_0^t\|u(\tau)\|_{B^s_{\infty,\infty}}^qA(\tau)^{2q}d\tau\Bigr)^\f1qB(t)^\f{q-1}q.\label{3.25}
\end{align}
Thus,  combining  \eqref{3.21}-\eqref{3.25} with (\ref{3.19}) and  using the Young's inequality lead to the result that
\begin{align}
A(t)^{2q}\le C(A(0)^{2q}+ \|u_0\|_2^{2q}+\|b_0\|_2^{2q})+C\int_0^t\|u(\tau)\|_{B^s_{\infty,\infty}}^qA(\tau)^{2q}d\tau.\nonumber
\end{align}
This together with  the Gronwall inequality  yields that
\begin{align}\label{3.26}
\sup_{t\in[0,T]}A(t)^{2q}\le C(A(0)^{2q}+ \|u_0\|_2^{2q}+\|b_0\|_2^{2q})\exp\Bigl(C\int_0^t\|u(\tau)\|_{B^s_{\infty,\infty}}^qd\tau\Bigr).
\end{align}
By means of (\ref{3.20}) and (\ref{3.26}), it follows
that there exists $\tilde{\rho}>\f12$ such that
$$
\sup_{t\in[0,T_0]}(\|u(t)\|_{H^{\tilde\rho}}+\|b(t)\|_{H^{\tilde\rho}})<+\infty,
$$
by Sobolev embedding  $B^{\rho}_{2,\infty}(\R^3)\hookrightarrow H^{\tilde\rho}(\R^3)$
with $\rho>\tilde{\rho}$ and $B^{1-\frac1q}_{2,\infty}(\R^3)\hookrightarrow H^{\tilde\rho}(\R^3)$
with $q>2$.
  Thus, the standard Picard's  method \cite{F-Kato, Kato}  ensures that the solution
$(u,b)$ can be extended after $t=T_0$ in the class ${\cal
X}(0,T_0)$. This completes the proof of Theorem 1.1.\endproof

\vspace{0.2cm}

\no\textbf{Acknowledgements}\quad  The authors thank the referees and the
associated editor for their invaluable comments and suggestions
which helped improve the paper greatly.  Q. Chen, C. Miao and Z.Zhang  were
supported by the NSF of China
under grant No.10701012, No.10725102 and No.10601002.

\end{document}